\documentclass{article}
\RequirePackage{etex}

\usepackage[margin=1in]{geometry}
\usepackage[english]{babel}
\usepackage[utf8]{inputenc}
\usepackage{subcaption}
\usepackage{amsmath}
\usepackage{amssymb}
\usepackage{amsfonts}
\usepackage{amsthm}
\usepackage{mathrsfs}
\usepackage[all]{xy}
\usepackage[pdftex]{graphicx}
\usepackage{color}
\usepackage{cite}
\usepackage{url}
\usepackage{indent first}
\usepackage[labelfont=bf,labelsep=period,justification=raggedright]{caption}
\usepackage[english]{babel}
\usepackage[utf8]{inputenc}
\usepackage{hyperref}
\usepackage[colorinlistoftodos]{todonotes}
\usepackage{tkz-fct}
\usepackage{tikz}
\usetikzlibrary{calc}
\usepackage{multicol}
\PassOptionsToPackage{dvipsnames,svgnames}{xcolor}
\usepackage{textcomp}

\DeclareMathOperator{\sgn}{\text{sgn}}
\DeclareMathOperator{\im}{\text{im}}
\DeclareMathOperator{\Span}{Span}
\DeclareMathOperator{\dis}{dis}

\def\multiset#1#2{\ensuremath{\left(\kern-.3em\left(\genfrac{}{}{0pt}{}{#1}{#2}\right)\kern-.3em\right)}}

\theoremstyle{plain}
\newtheorem{thm}{Theorem}

\newtheorem{lemma}[thm]{Lemma}

\theoremstyle{definition}

\theoremstyle{remark}
\newtheorem{rem}[thm]{Remark}

\numberwithin{equation}{section}
\numberwithin{thm}{section}

\usepackage{arxiv}

\usepackage[utf8]{inputenc} 
\usepackage[T1]{fontenc}    
\usepackage{hyperref}       
\usepackage{url}            
\usepackage{booktabs}       
\usepackage{amsfonts}       
\usepackage{nicefrac}       
\usepackage{microtype}      
\usepackage{lipsum}

\title{Irreducibility of Markov Chains on simplicial complexes, the  Spectrum of the Discrete Hodge Laplacian and Homology}

\author{
Marzieh Eidi\\
 Center for Scalable Data Analytics and Artificial Intelligence, Leipzig University\\
 Max Planck Institute for Mathematics in the Sciences\\
  \texttt{meidi@mis.mpg.de} 
   \AND
Sayan Mukherjee\\
 Center for Scalable Data Analytics and Artificial Intelligence,  Leipzig University\\
 Max Planck Institute for Mathematics in the Sciences\\
 Duke University\\
 \texttt{sayan.mukherjee@mis.mpg.de} 
}

\begin{document}

\maketitle

\begin{abstract}
Random walks on graphs are a fundamental concept in graph theory and play a crucial role in solving a wide range of theoretical and applied problems in discrete math, probability, theoretical computer science, network science, and machine learning. The connection between Markov chains on graphs and their geometric and topological structures   is the main reason why such a wide range of theoretical and practical applications exist. Graph connectedness ensures irreducibility of a Markov chain. The convergence rate to the stationary distribution is determined by the spectrum of the graph Laplacian which is associated with lower bounds on graph curvature. Furthermore, walks on graphs
are used to infer structural properties of underlying manifolds in data analysis and manifold learning.  However, an important question remains: can similar connections be established between Markov chains on simplicial complexes and the topology, geometry, and spectral properties of complexes? Additionally, can we gain topological, geometric, or analytic information about a manifold by defining appropriate Markov chains on its triangulations? These questions are not only theoretically important but answers to them provide powerful tools for the analysis of complex networks that  go beyond the analysis of pairwise interactions. In this paper, we provide an integrated overview of the existing results on random walks on simplicial complexes, using the novel perspective of signed graphs. This perspective sheds light on previously unknown aspects such as irreducibility conditions. We show that while up-walks on higher dimensional simplexes can never be irreducible, the down walks become irreducible if and only if the complex is orientable. We believe that this new integrated perspective can be extended beyond discrete structures and enables exploration of classical problems for triangulable manifolds.


\end{abstract}

\keywords{Signed Graphs, Random walks, Normalized Laplacian, Simplicial complexes, Laplacian spectrum, (Co)Homology, Orientability, Disorientability, Stationary Distribution, Perron-Frobenius Theorem}

\section{Introduction}

Stochastic processes on mathematical structures can reflect their topological, geometric and spectral properties. This ranges from Brownian motion on smooth Riemannian manifolds \cite{Hsu} to random walks on graphs. Random walks on graphs and Riemannian manifolds exhibit similar properties
and by exploring the convergence properties of random walks, one can investigate the approximation of Riemannian manifolds by graphs \cite{parva}. The discrete generators of  random walks on graphs approximate the Laplace-Beltrami operator on the underlying manifold and has been used 
in developing algorithms that learn
the underlying manifold of a data set \cite{belkin,LLE,isomap,tsne}. For graphs the limiting behaviour of the process tells us about connectedness (zero dimensional homology), and the rate of convergence is related to the  spectral gap of the Laplacian \cite{Lov}. 
Random walks on graphs have been used to prove the Cheeger inequality for graphs \cite{Chung2}. The discrete Cheeger inequality is one of the most fundamental quantities in discrete optimization, spectral
graph theory and the analysis of Markov chains. It relates the spectral gap of the normalized graph Laplacian to the
the isoperimetric constant. This inequality for graphs is the discrete version of Cheeger's inequality for compact Riemannian manifolds which was proved by Jeff Cheeger in 1970 \cite{cheeg1}. The isoperimetric constant
is a positive real number defined in terms of the minimal area of a hypersurface that divides the manifold into two disjoint pieces. This inequality relates the first nontrivial eigenvalue of the Laplace–Beltrami operator on the manifold to its Cheeger constant.  The discrete Cheeger inequality has  many applications in graph clustering \cite{cheeg3}, 
expander graphs\cite{cheeg2}, analysis of Markov chains, and image segmentation \cite{cheeg1}. 

Furthermore, random walks on graphs are widely utilized for inference of geometric properties of the graph. In mathematics a striking example is the generalized Ricci curvature proposed by Ollivier \cite{Ollivier}; he defined a notion of Ricci curvature for Markov chains  on metric measure spaces. This ranges from Rimannian manifolds when considering Riemannian volume measure to discrete structures such as graphs. Then  Lin-Lu-Yau \cite{Yau} developed this curvature notion in graphs and showed that  a  lower Ricci curvature bound  implies various geometric and analytic properties; the first non-zero eigenvalue of the Laplacian, which determines the rate of convergence, is bounded from below by the positive lower bound of this curvature. Interestingly, all these results originated from corresponding known theorems in geometric analysis of Riemannian manifolds \cite{jbook}. Due to these strong implications, in the applied side, Ollivier Ricci curvature notion has been among the most popular methods used for the analysis of complex networks ranging from graphs to directed hypergraphs in the past few years \cite{Faro, market, eidi}.

Graphs are  one dimensional simplicial complexes where the state space of the Markov chain consists of the graph vertices, allowing for random movement between any two vertex as long as they are connected by an edge.
To extend graph random walks to higher dimensional simplicial complexes,  different ideas have been proposed where their asymptotic behaviours are related to the higher dimensional spectral gap of the Hodge Laplacian and can reflect the non-trivial high-dimensional homology in the corresponding dimensions
\cite{Ros, Sayan,Schaub_2020}. However, there are two main differences between random walks on graphs and simplicial complexes in arbitrary dimensions. The first difference is that in a graph moving on the vertices is based on the  upper adjacency connections between them; namely the walk is  via the common edges between these vertices. However when we start walking from the $d$-dimensional simplexes to reach other $d$-dimensional simplexes, for $1\leq d < N$ we can move both upward and downward in dimension; namely through $(d+1)$-simplexes or $(d-1)$-simplexes. Also for $d=N$, we can just move downward through $(N-1)$-simplexes. 
The second difference which has been more difficult to deal with is the matter of orientation; for every non-zero dimensional simplex, we have two possible orientations; an orientation for a $d$-simplex ($d \succ 0$) is an equivalence class of
orderings of its vertices, where two orderings are equivalent if they differ by an even
permutation. Therefore based on the ideas in the literature for both up and down (and full) random walks \cite{Ros, Sayan,Schaub_2020}, the state space of Markov chain is no longer the set of $d$-simplexes, but the "oriented" $d$-simplexes which doubles the size of the number of $d$-simplexes. The orientation complicates the connection of the limiting behaviour of random walks to the $d$-dimensional spectral gap of the Laplacian and the existence of non-trivial homology.  The main challenge is that in higher dimensions we deal with the limit of the difference between two probability distributions and not a single one as in the graph case. This 
difference is much harder to analyse and 
in this setting the limiting behaviour of Markov chains is no longer independent of the starting distribution. There is no clear understanding of when and if  Markov chains on simplicial complexes are irreducible. Furthermore, to prevent periodicity, similar to the (bipartite) graph case, we need to use a specific laziness parameter on the random walks; for the down-walk (i.e. via shared faces) this parameter is dependent on the degree of the simplexes \cite{Sayan} and  for the up walk it depends on the dimension of the simplicial complex\cite{Ros}. Recall that both irreducibility and aperiodicity are needed for a finite state Markov chain to have a unique stationary distribution (i.e. to be ergodic). This stationary distribution then is used to connect to the geometry and (co)-homology of the underlying structure. In terms of applications,
random walks on graphs and simplicial complexes are the main tools in tackling diverse complex real-words problems from ranking  web-pages in Google page rank algorithm \cite{Page1,Page2} to signal processing and flow network decomposition\cite{Schaub_2020}. 
Considering the numerous theoretical implications and practical applications that random walks on
graphs and simplicial complexes have, it is very desirable to develop an understanding of the analogous properties on graphs and in particular irreducibility of the walks on higher dimensions in simplicial complexes. 

There are several reasons why we care about the irreducibility of finite Markov chains \cite{Markov1, Markov2}:
\begin{itemize}
  
  \item Predictability: An irreducible Markov chain ensures that every state can be reached from any other state. This property ensures that the chain does not get stuck in certain states and eventually converges to its stationary distribution and  is predictable.  There are no isolated or disconnected states and allows for accurate predictions about the future states of the system based on the current state.
\item Ergodicity: Irreducibility is a key requirement for a Markov chain to be ergodic. An ergodic Markov chain has the property that the long-run average of any observable variable converges to a fixed value. This property is important for studying the stability and equilibrium of systems modeled by Markov chains and we can use it to understand how the system will behave in the long run.

\item Applications in various fields: Finite state Markov chains are widely used in mathematical fields such as Rimannian geometry, probability theory, statistics and graph theory. They are also foundational for algorithms used in 
Natural Language Processing (NLP), physics, finance and economics, biology,
optimization and algorithms in applied
computer science. Irreducible Markov chains have been used in learning algorithms, which aim to uncover the underlying structure of high-dimensional data lying on a Riemannian manifold. By constructing Markov chains that capture the local relationships between data points, one can perform tasks such as dimensionality reduction and clustering on a manifold. Markov chain Monte Carlo (MCMC) methods  which are based on
 irreducible Markov chains are widely used for sampling from complex distributions and solving optimization problems on Riemannian manifolds \cite{mcmc,mcmc2, mcmc1}.

\end{itemize}

Simplicial complexes have become popular models for representing geometric, topological, and higher-order interactions:
\begin{itemize}
    \item 

They provide a more flexible and expressive representation than graphs and unlike graphs have a natural hierarchical structure. While graphs are limited to representing pairwise relationships between vertices, simplicial complexes can represent higher-order relationships between groups of vertices. This allows for the modeling of more complex structures and interactions and  their hierarchical structure allows for the exploration of multi-scale features and phenomena.

\item Simplicial complexes enable the study of topological properties and structures. Graphs lack the notion of boundaries and higher-dimensional geometric properties, which are crucial in many fields such as topology and geometry.

\item Simplicial complexes come with efficient algorithms and computational tools. Many techniques and algorithms, such as homology, persistent homology in  topological data analysis, and simplicial complex-based machine learning, have been developed in the past 10 years for simplicial complexes. This has made them more and more suitable for analyzing and interpreting complex datasets and networks.

\end{itemize}

Combining irreducible Markov chains and simplicial complexes is the main aim of this article to address the unclear aspects about irreducibility of Markov chains on simplicial complexes by assembling various puzzle pieces from the literature to form an overall picture that enables us to clarify conditions that ensure irreducibility. Our main tool is signed graphs and the known connection between their spectrum with the spectrum of the Hodge up/down Laplacian on simplicial complexes. In the next section we state the definitions used throughout the paper.

\section{Definitions and preliminaries}

\textit{Simplicial complexes, boundary, and co-boundary maps:}
An abstract simplicial complex $K$ on a finite set $V$, $V=\lbrace1,...n\rbrace$ is a collection of subsets
of $V$ , which is closed under inclusion. A $d$-simplex of $K$ is an
element of cardinality $d+1$. $0$-simplexes are usually called vertices and $1$-simplexes
edges. The collection of all $d$-simplexes of simplicial complex $K$ is denoted by
$S_d$. The dimension of a $d$-simplex is $d$, and the dimension of a complex $K$
is the maximum dimension of a simplex in $K$. In this article we assume the dimension of $K$ is $N$. The simplexes which are maximal
under inclusion are called facets. Let  $[S_d]$ be the set of all oriented $d$-simplexes. For $1 \leq d \leq N$, there are two opposite orientations  of a d-simplex $\sigma$  by $[\sigma]$ and $-[\sigma]$ (note that both are in $[S_d]$).  

The $d$-th chain group $C_d(K)$ of $K$ is a vector space with the basis $S_d$. We consider the setting with real coefficients. The {\em boundary map } 
$\partial_d: C_d(K) \rightarrow C_{d-1}(K)$ is a linear operator defined by 
\begin{equation} 
\partial_d[i_0,...i_d]=\sum_{j=0}^{d} (-1)^j [i_0,. i_{j-1},i_{j+1}..i_d].
\end{equation}

The sequence ($C_d (K), \partial_d$) is the  chain complex of K, meaning that $\partial_d \partial_{d+1} =0 $ for all
$d$. Therefore we can consider the kernel and image of the boundary maps, denoted by $\ker \partial_d$  and  $\im \partial_{d+1}$.
The kernel is called the space of cycles and the image is called the space of boundaries. We denote $H_d= \frac{\ker \partial_d}{\im \partial_{d+1}}$ as the (real) $d$-homology.
 Every boundary is also a cycle and non-trivial cycles are those which are not boundaries and correspond to the non-trivial homology.

The $d$-th co-chain group $C^d(K)$ is defined as the dual of the chain group $C_d(K)$ with real coefficients. These are the functions that satisfy $f(-[\sigma])=-f([\sigma])$ for every oriented simplex $[\sigma]$.

The coboundary map $\delta_d: C^d(K)\rightarrow C^{d+1}(K)$ can be defined by the following linear operator

\begin{equation}
\delta_df([i_0,...i_d])=\sum_{j=0}^{d+1} (-1)^j f[i_0,. i_{j-1},i_{j+1}..i_d] 
\end{equation}
for any $f \in C^d(K)$. 

The  co-chain complex of K is ($C^d (K), \delta_d$), meaning that $\delta_{d} \delta_{d-1} =0 $ for all
$d$. The kernel and image of the coboundary map are denoted $\ker \delta_d$ and $\im \delta_{d-1}$.
The kernel is the space of co-cycles (closed forms) and the image is the space of co-boundaries (exact forms). The $d$-th cohomology
is $H^d= \frac{\ker \delta_d}{\im \delta_{d-1}}$.
Closed forms always contain exact forms and existence of non-trivial closed form, i.e those which are not exact corresponds to the non-trivial cohomology. $\delta_d$ can be viewed as the dual of the boundary map $\partial_{d+1}$.  For more details on simplicial homology and cohomology the reader is
referred to \cite{hatcher}. \\ 
We can choose a positive definite inner product $\langle \cdot, \cdot \rangle_d$ on $C^d$. There is a one-to-one correspondence between weight functions on the set of simplexes
and possible scalar products on co-chain groups  such that elementary co-chains are orthogonal \cite{Horak}. Therefore the inner product is given by a weight function $w: K\rightarrow (0, \infty)$ with
\begin{equation}
    \langle f, g \rangle = \sum_{\sigma \in S_d} w(\sigma) f(\sigma) g(\sigma)
\end{equation}

For finite dimensional simplicial complexes we can then define the adjoint $\delta^*_d: C^{d+1}(K)\rightarrow C^d(K)$ of the coboundary operators $\delta_d$ as:
\begin{equation}
    (\delta_d f_1, f_2)_{C^{d+1}}=(f_1, \delta^*_df_2)_{C^{d}}
\end{equation}
for every $f_1 \in {C^{d}(K)}$ and $f_2 \in {C^{d+1}(K)}.$

If $w$ is the weight function, the degree of a $d$-simplex $\sigma$ is equal to the sum of the
weights of all simplexes that contain $\sigma$ in their boundary
\begin{equation}
    \deg \sigma := \sum_{\tau_{d+1}: \sigma \in \partial \tau } w(\tau)
\end{equation}

\textit{Discrete Hodge Laplacians:}
Three discrete operators will be used in this paper. 

For further details of discrete Laplacians and the connection between continuous and discrete Laplacians see \cite{Horak,Max}.

 Let the weight function be the constant $1$ on the set of all of the simplexes, for each $d$-simplex, $ 0 \leq d \leq N$. We can define the following three operators on $C^d(K)$.
\begin{itemize}
    \item  The combinatorial $d$-th up-Laplace operator $L_d ^{up}:= \delta^*_{d} \delta_d$.
    \item The combinatorial $d$-th down-Laplace operator $L_d ^{down}:= \delta_{d-1} \delta^*_{d-1}$.
    \item The combinatorial (full) $d$-Laplacian  $L_d= L_d ^{up} +L_d ^{down}$.
\end{itemize}
All of these operators are self-adjoint and non-negative and therefore have non-negative real eigenvalues. For $d=N$, $L_d ^{up}$ is zero and for $d=0$, $L_d ^{down}$ is zero.

From the above definitions of Laplace operators it holds that $L_d ^{up}$ is uniquely determined by its restriction on the $(d+1)$-skeleton of
$K$ and $L_d ^{down}$ is determined by its $d$-skeleton. Therefore, it suffices to observe pure $(d+1)$- simplicial complexes for the up Laplacian and pure $d$-simplicial complexes for the down Laplacian. 
Recall that a simplicial complex is called 
pure if all facets have the same dimension. \\
This is a fundamental fact that will be used in this paper when dealing with Laplace operators for arbitrary $d$.

We can extend Laplacian definitions to arbitrary positive weight functions and for the weighted up Laplacian we obtain \cite{Horak}:
\begin{equation}
(L_d^{up}f)[\sigma]:= \sum_{\tau_{d+1: \sigma \in \partial \tau }} \frac{w(\tau)}{w(\sigma)} f([\sigma])+\sum_{\sigma'\neq \sigma, \exists \tau_{d+1}:\sigma, \sigma' \in \partial \tau\ }\frac{w(\tau)}{w(\sigma)} \sgn ([\sigma], \partial [\tau]) \sgn  ([\sigma'], \partial[ \tau])f([\sigma']).
\end{equation}
Here $\partial[\tau]$ is the induced orientation of $\tau$ on its faces (and in particular on $\sigma$ and $\sigma'$) by the boundary operator $\partial$. Note $\sgn ([\sigma], \partial[ \tau])$ is positive if the assigned orientation to $\sigma$ is the same as the induced orientation from $\tau$ by the boundary operator on $\sigma$ and is negative otherwise.
A similar formula can be obtained for the weighted down Laplacian 
 \begin{equation}
(L_d^{down}f)[\sigma]:= \sum_{\rho_{d-1: \rho \in \partial \sigma }} \frac{w(\sigma)}{w(\rho)} f([\sigma])+\sum_{\sigma', \sigma \cap \sigma'= \rho }\frac{w(\sigma')}{w(\rho)} \sgn ([\rho], \partial[ \sigma]) \sgn([\rho], \partial[ \sigma'])f([\sigma']).
\end{equation} 
Note that the combinatorial up/down Laplacians can be obtained by these formulas when considering the constant $1$ weight function on the set of simplexes.   
A special case of the above formulas are the "normalized Laplacians" where for each simplex $\sigma$ such that $\sigma$ is not a facet, $\deg \sigma= w(\sigma)$ and the weights of all the facets are one.  We use $\Delta$ to denote the 'normalized' Laplacian. The normalized Laplacian is suitable for problems related to random walks, similar to the graph case, as well as the
expanders. For more details on the advantages of this weight function see \cite{Horak}.

\textit{The spectrum of Laplacians and Topology:}
A fundamental theorem which connects the spectrum of Laplacian to the topology of the simplicial complexes is Eckmann's theorem, also known as the discrete Hodge theorem \cite{Eckmann}. This theorem states an isomorphism between the kernel of the Laplacian (harmonic functions)  with the cohomolgy of the simplicial complex. The multiplicity of the eigenvalue zero of $L_d$ is equal to the Betti number, the dimension of the $d$-th  cohomology. For finite dimensional simplicial complexes and field coefficients the $d$-th  cohomology and $d$-th  homology are isomorphic.
 The following important decomposition is a main result of this theorem:
\begin{equation}
    C^d= \im \delta^*_{d} \oplus 
 H^d \oplus \im \delta_{d-1}
\end{equation}

We note that although the kernel of the Laplacian depends on the choice of scalar product on the space of forms (i.e. depends on the weight function defined on the simplexes), Eckmann's theorem does not. All  three Laplacians decompose with respect to
this decomposition. By definition, $\ker L_d= \ker L_d^{up} \bigcap \ker L_d^{down}$. It holds that: 
$\ker L_d^{up}=\ker \delta_{d}$ and $\im L_d^{down}=\im \delta_{d-1}$. The spectrum of $L_d$ depends on the scalar product and the non-zero spectrum of $L_d$ is the union of the non-zero spectrum of $L_d^{up}$ and $L_d^{down}$. Also, there is a one-to-one correspondence between the non-zero spectrum of $L_d^{up}$ with the non-zero spectrum of $L_{d+1}^{down}$. This is why we usually consider just one of them when dealing with the non-zero spectrum (usually the up-Laplacian is considered). 

 \textit{Connected simplicial complex:} We will work with the up-adjacency matrix for the up-Laplacian and the down-adjacency matrix for the down-Laplacian. This will require two definitions of connectedness of a simplicial complex.
\begin{itemize}
\item  Up-connectedness: An $N$-complex $K$ is called up $d$-connected $(0 \leq d \leq N-1)$ if for every two $d$-simplexes $\sigma$, $\sigma'$ there exists a chain of two by two upper adjacent $d$-simplexes that connect $\sigma$ to $\sigma'$. Such chains define equivalence classes of $d$-simplexes, and the
simplexes in the same class are in the same up $d$-connected component.

\item  Down-connectedness: An $N$-complex $K$ is called down $d$-connected $(1 \leq d \leq N)$ if for every two $d$-simplexes $\sigma$, $\sigma'$ there exists a chain of two by two lower adjacent $d$-simplexes that connect $\sigma$ to $\sigma'$. The equivalence classes are the down $d$-connected components.  
\end{itemize}


\textit{Disorientability and Orientability of simplicial complexes:} A disorientation of an $N$-complex $K$ is a choice of orientation of its $N$-simplexes, that whenever two arbitrary $N$-simplexes intersect in a $(N-1)$-simplex, they induce the same orientation on the $(N-1)$-simplex. If $K$ has a disorientation it is said to be disorientable. This notion was originally introduced as the higher dimensional analogue of bipartiteness \cite{Ros}.  We note that disorientability is defined for the maximum dimension of the complex \cite{Sayan}. Bipartite graphs are the simplest examples of  disorintable graphs and in fact a graph is disorientable if and only if it is bipartite. 

A $N$-complex is orientable  if there is
a choice of orientations of its $N$-simplexes, so that simplexes intersecting in a codimension one simplex
induce opposite orientations on it.

Note that disorientability and orientability are not the opposite of each other and a simplicial complex can be both orientable and disorientable but in an orientable simplicial complex, the degree of $(N-1)$- simplexes is at most $2$ while for disorientable ones we have no such restriction.  

\section{Random Walks on graphs and higher dimensional Simplicial complexes} \label{rwsc} 
In this section we  state existing constructions of random walks and their properties. We start with classical results for graphs, $1$-d simplicial complexes and then we go to the higher dimensional case.  In the second subsection we show how the relation between the Laplacian on signed graphs and the Hodge Laplacian coupled with irreducibility conditions allow us to get a simplified picture of the existing results. This enables us to answer the following important question:

Is it possible to define a graph-type random walk on $d$-simplexes  that can reflect non-trivial $d$-th homology? By graph-like we mean a Markov chain with states that are $d$-simplexes (not oriented $d$-simplexes). In graphs vertices  are the states and  do not have orientations.   
\paragraph{Random Walks on Graphs:}
The state space of the Markov chain $M$ consists of all the vertices of the graph $G$. To prevent periodicity, we usually work with lazy random walks where we have a positive chance $\alpha$ of staying at each state.  Recall that for a Markov chain, aperiodicity is a necessary condition  to have a  unique stationary distribution. The  transition between the states is defined as following: the walker starts
at a vertex $v_0$, and at each step remains in place with probability $\alpha$
or moves to each of its 
neighbors with probability $\frac{1-\alpha}{\deg(v_0)}$. If $p_t^{v_0}(v)$ is the probability of finding the walker at $v$ after t-step walk starting from $v_0$, then the following results are classic:
\begin{itemize}
\item The (normalized) Laplacian $\Delta_0$ is the discrete generator of Markov chain: \begin{equation}\label{graph}
    M= I- (1-\alpha)\Delta_0  
\end{equation}
\item If the graph $G$ is finite, then $p^{v_0}_\infty = lim_ {t \rightarrow\infty} p^{v_0}_t$
exists and is independent of the starting distribution if and only if graph is connected.
\item  The rate of convergence to the stationary distribution, which is unique due to aperiodicity, is given by 
\begin{equation}
\dis(p_t, p_{\infty}  )= O(1-\alpha \lambda(G))^t. 
\end{equation}
where $\lambda(G)$ is the spectral gap of G, (i.e. the first non-zero eigenvalue of the graph Laplacian).
\end{itemize}
 
Coming back to the previous discussion of random walks and the Laplacian  of simplicial complexes, we note that except for the minimum and maximum dimensions of a simplicial complex the full Laplacian can be decomposed to the up-Laplacian and down-Laplacian. Therefore for the middle dimensions we can define both up and down walks while for the minimum dimension we can just define the up walks (such as known walks on the vertices of a graph) and in the maximum dimension we can just define the down walk.
 But what is common in all of these approaches \cite{Ros, Sayan, Schaub_2020} is that, for the random walks beyond vertices, the state space of the Markov chains are "oriented" simplexes  which has the twice number of elements of the simplexes in that dimension.

\subsection{Existing approaches for Random Walks on Simplicial complexes and their challenges}

\paragraph{Up-random walk}
Rosenthal and Parzanchevski \cite{Ros} defined Neumann random walks on the oriented $(N-1)$- dimensional simplexes of a $N$-dimensional simplicial complex, $[S_{N-1}]$; this state space has twice the $\#S_{N-1}$ elements as each simplex of dimension higher than zero has two orientations.

\begin{itemize}
\item [(1)] Co-neighbors:  two oriented $(N-1)$-simplexes are called co-neighbors, which we denote as $[\sigma] \uparrow [\sigma'] $
if they share a $N$- coface and $[\sigma]$ and $[\sigma']$ induce the same orientation on their shared face.
\item[(2)] Transition matrix $P$: the transition matrix is the time-homogeneous Markov chain $M$ on the state space $[S_{N-1}]$ with transition probabilities
\[ P_{[\sigma],[\sigma']} = \text{Prob}([\sigma] \to [\sigma'])  = \begin{cases}
 p & \mbox{when }  [\sigma'] = [\sigma]  \\
 \frac{1-p}{N \deg(\sigma)} & \mbox{when }  [\sigma'] \uparrow [\sigma]  \\
0 & \mbox{otherwise} 
\end{cases} \]
\end{itemize}
where deg$(\sigma)$ is the number of its N-cofaces.\\
One randomly moves between upper adjacent oriented $(N-1)$-simplexes as long as they induce the same orientation on their common $(N-2)$-face. 
For the transition  operator on a simple example, consider the one depicted in Fig.1. We start the random walk from the  oriented edge $[e_1]$ and with probability $1/2$ we stay at $[e_1]$ or based on the neighboring condition, we can jump to one of the oriented edges: $[e_2]$,  $[e_3]$,  $[e_4]$,  $[e_5]$ each with equal probability$=1/8$. Note that all of $[e_1], [e_2],[e_4]$ have the same starting vertex and $[e_1], [e_3],[e_5]$ have the same ending vertex and this guarantees the neighbouring condition  among these oriented edges. 
\begin{figure}[ht]
\centering
\includegraphics[width=0.35\textwidth]{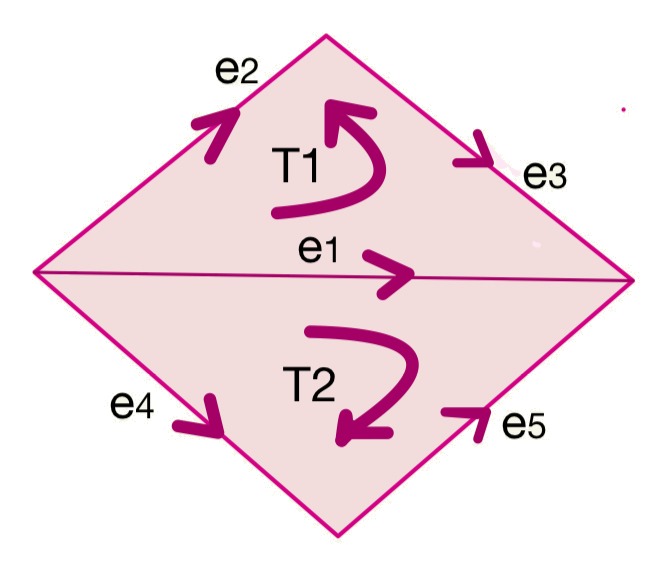}
\caption{Up random walk (i.e. via triangles) on "oriented" edges of a 2-d simplicial complex. 
}
\end{figure}

 A key idea in this work was that if $p_t^{[e_0]}([e])$ is the probability of finding the walker at [e] after $t$-step walk starting from $[e_0]$, the limit of "the normalized expectation process", is connected to the spectrum of the up-Laplacian and can reflect the non-trivial $1$-homology. The normalized expectation process is defined as
$$\tilde{\mathcal{E}}_t := \left(\frac{4}{3}\right)^t \left[(p_t^{[e_0]}([e])-p_t^{[e_0]}(-[e]))\right]$$ $\tilde{\mathcal{E}}_t$ is seen as a co-chain in $C^1$ and what matters when it comes to topology is to determine when this co-chain becomes a co-cycle and co-boundary.  
Also the reason for using the scaling factor $(\frac{4}{3})^t$, is that in the limit, for any starting $[e_0]$, $p_t^{[e_0]}([e])$ and $p_t^{[e_0]}(-[e])$  (i. e the probabilities of reaching $[e]$ and its reverse) become arbitrarily close and therefore their difference becomes zero.  Also in this walk specific thresholds for laziness, which depends on the dimension of the complex, is needed to be sure that the limiting behaviour can detect homology. Their results are not only obtained for the edges in $2$-complexes but for $(N-1)$-simplexes of a $N$-complex, for arbitrary $N$ and arbitrary laziness parameter $p$ . For this general case, the $\frac{4}{3}$in the above definition is substituted by $\frac{N}{p(N-1)+1}$.  
Furthermore, the evolution of the normalized expectation process over time is given by$\tilde{\mathcal{E}}^{[\sigma_0]}_{t+1}= A \tilde{\mathcal{E}}^{[\sigma_0]}_t$, where $A$ is the transition operator acting on  $C^{N-1}$ which is the transpose of $M$ with respect to a natural choice of basis for $C^{N-1}$. 

\begin{thm}\label{Ros} 
With the above assumptions the following are the main results \cite{Ros}: 
\begin{itemize}
\item[i.] $$A=\frac{p(N-1)+1}{N}I- \frac{1-p}{N} L_{N-1}^{up}$$ 
where $L$ is the Laplacian  with  a wight function that for each $(N-1)$- simplex $\sigma$, $w(\sigma)=\frac{1}{deg(\sigma)}$ and $"1"$ elsewhere. 
\item[ii.]  If $\frac{N-1}{3N-1}<p<1 $, then 
$$\tilde{\mathcal{E}}_{\infty}^{[\sigma_0]}= lim_{t\rightarrow \infty} \tilde{\mathcal{E}}_t^{[\sigma_0]}$$
is an exact co-chain for every starting oriented simplex $[\sigma_0]$ if and only if $H_{N-1}(K)=0$ and  the dimension of $H_{N-1}(K)$ equals the dimension of 
$$\Span \{\mathbb{P}_{\ker \delta^*_{N-2}}(\tilde{\mathcal{E}}_{\infty}^ {[\sigma]})\mid \sigma \in S_{N-1}\} $$
where $\mathbb{P}$ denotes the projection map onto the $\ker \delta^*_{N-2}$. 
 
\item[iii.] If  $p\geq 1/2 $, then the rate of convergence is controlled
by the spectral gap $\lambda_{N-1}$ of $L_{N-1}^{up}$:
\begin{equation}
\dis(\tilde{\mathcal{E}}_t^{[\sigma_0]}, \tilde{\mathcal{E}}_{\infty}^{[\sigma_0]})= O\left((1-\frac{1-p}{p(N-1)+1} \lambda_{N-1}\right)^t. 
\end{equation}
 Also when $p= \frac{N-1}{3N-1}$, $\tilde{\mathcal{E}}_{\infty}$ is exact for every starting $[\sigma_0]$ if and only if the $(N-1)$-homology is trivial and in addition there is no disorientable $(N-1)$-components.\\
  More generally, for $\frac{N-1}{3N-1}<p<1 $, $\tilde{\mathcal{E}}_{\infty}^{[\sigma_0]}$ is closed and the same holds  when $p= \frac{N-1}{3N-1}$ unless $K$ has a disorientable $(N-1)$ component.  
\end{itemize}
\end{thm}

\begin{rem}
We note that since this walk is defined on the ($N-1$)- oriented simplexes, and this is neither the minimum (we exclude the graph case) nor the maximum dimension, the zeroes of $L^{up}_{N-1}$ are not necessarily the zeroes of the full Laplacian and therefore are not directly related to the $(N-1)$-homology. In \cite{Ros} for generalising the spectral gap of the normalized graph Laplacian, the authors introduce the followings:  

The spectral gap $\lambda(K)$ and essential gap $\Bar\lambda(K)$ of $L_{N-1}^{up}$ in a finite $N$-dimensional complex $K$ are
\begin{equation} 
\lambda(K) = \min \mbox{Spec} \left( L_{N-1}^{up} \mid_{\ker L^{down}_{N-1}} 
\right) =
 \min \mbox{Spec} \left(L_{N-1}\mid_{ \ker L_{N-1}^{down}} \right)
\end{equation}
\begin{equation} 
\overline{\lambda}(K) = \min \mbox{Spec} \left(L_{N-1}^{up}\mid_{ \im L_{N-1}^{up}} \right)= \min \text{Spec}
\left (L_{N-1}\mid_{\im L_{N-1}^{up}} \right)
\end{equation}
 In these definitions, $\lambda = 0$ whenever the $(N-1)$-th homology is non-trivial and $\overline{\lambda}$ does not vanish. This has additional advantages when $(N-1)$-th homology is trivial. Therefore in the above theorem, the convergence rate to the stationary distribution, depending whether it is exact or not, is determined by either of spectral or essential gaps of the Laplacian.   

\end{rem}

Although the above theorem is similar to the one already mentioned for graphs, analysis of this process and understanding scaling limits is challenging because in expectation the process is centered at zero. More importantly, the limiting behavior of the process is dependent on the staring point; namely, as mentioned in \cite{Ros}, there are complexes (even $(N-1)$-connected ones) with nontrivial $(N-1)$-
homology, such that the $\tilde{\mathcal{E}}_{\infty}^{[\sigma_0]}$ is exact for a carefully chosen $[\sigma_0]$. This dependency prevents having a clear idea about the irreducibilty of the walk and recognition of the necessary and/or sufficient conditions to achieve that.

These two difficulties are the main difference between these walks and well understood walks on  graphs, where we consider the evolution of a single probability measure of Markov chains on connected graphs that are irreducible. 
\begin{rem}
    In \cite{Ros} the authors use the transition operator $A$, to derive the above theorem. However since $A$ and $M$, the Markov chain given by the transition probabilities described above, are transpose of each other (with respect to a natural choice of basis for $C^{N-1}$), the same results can be obtained by substituting $A$ in the above theorem with $M$, and that  corresponds to the normalized up-Laplacian on the simplicial complex.    
\end{rem}

\paragraph{Down random walk}
In \cite{Sayan} a Dirichlet walk was defined that is related to the spectrum of the down $d$-Laplacian; here the state space is all of the oriented 
$d$-simplexes (namely the set $[S_d]$)  plus a "death-state" $\Theta$.

\begin{itemize}
\item[(1)] Neighbors:  two oriented $d$-simplexes $[\sigma]$ and $ [\sigma'] $ are called neighbors, which we denote as 
$[\sigma] \downarrow [\sigma'] $ if they share a $(d-1)$-face and induce opposite orientation on it.
\item[(2)] Transition matrix $P$: the transition matrix is the time-homogeneous Markov chain on the state space $[S_d] \cup \Theta$ with transition probabilities
\[ P_{[\sigma],[\sigma']} = \mbox{Prob}([\sigma] \to [\sigma'])  = \begin{cases}
p  & \mbox{when } \sigma \neq  \Theta \mbox{ and } [\sigma] = [\sigma'] \\
\frac{1-p}{(M-1)(d+1)}  & \mbox{when } \sigma \neq \Theta \mbox{ and }   [\sigma] \downarrow [\sigma'] \\

0  & 
\mbox{otherwise },
\end{cases} \]
where $M$ is the maximum degree of  $(d-1)$-simplexes in the complex.

Also the transition probabilities for switching the orientations are :
\[ P_{[\sigma],-[\sigma']} = \text{Prob}([\sigma] \to -[\sigma'])  = \begin{cases}
\frac{1-p}{(M-1)(d+1)}  & \mbox{when } \sigma \neq \Theta \mbox{ and }   [\sigma] \downarrow -[\sigma'] \\
0  & 
\mbox{otherwise },
\end{cases} \]

The probability of transitioning into a death state is equal to: 
$$1- \sum_{[\sigma] \downarrow [\sigma'], [\sigma] \neq [\sigma']} \mbox{Prob}([\sigma] \to [\sigma'])- \sum_{[\sigma] \downarrow -[\sigma']} \mbox{Prob}([\sigma] \to -[\sigma'])-p.$$
The probability of going from $\Theta$ to $\Theta$ is one. The transition probabilities are constrained to ensure $P$ is a stochastic matrix.
\end{itemize}

Consider the following example  of randomly jumping between oriented edges via their shared vertices. Since we want to examine the  down $1$-Laplacian, we do not care about the existence of the $2$-dimensional triangles $T_1, T_2$. 

\begin{figure}[ht]
\centering
\includegraphics[width=0.35\textwidth]{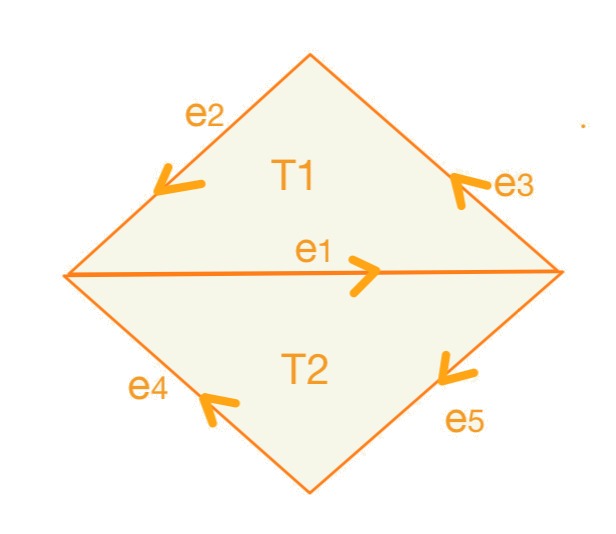}
\caption{Down random walk (i.e. via vertices) on "oriented" edges of a 2-d simplicial complex. 
}
\end{figure}
Starting at $[e_1]$, the walk has
probability $p$ of not moving and for each of the neighbors of $[e_1]$ the walk has probability $\frac{1-p}{2(M-1)}$ of jumping to that neighbor, where $M$ is the maximum degree of the vertices (=3 in this graph). If the number of neighbors of $[e_1]$ is less than $2(M-1)$,
then the sum of these probabilities is less than $1$. In this case, we
interpret the difference as the probability that the walker dies (i.e.,
the walker jumps to a death state from which it can never return).
 In order to connect this walk to the spectrum of the down $1$-Laplacian (and therefore to the $1$-th homology), the evolution of  
 $T P^t \nu$ was examined, which is the image of the marginal distribution under a
linear transformation $T$ ($\nu$ is the initial distribution) such that  for any function
$f:[S_1] \cup \Theta\rightarrow \mathbb{R}$, 
$$Tf[e]= f([e])- f(-[e])$$ 
(i. e.  $T$ enforces antisymmetry between the opposite orientations of an edge).  The spectrum, and in particular the spectral gap, of the down $1$-Laplacian and the limit of $CTP^t$ are related where $C$ is a normalizing constant, depending on $M$, that ensures $CTP^t$ has nontrivial limiting behavior. To this end the following propagation matrix was used
\[ B_{[\sigma'],[\sigma']} =  = \begin{cases}
p  & \mbox{when } [\sigma] = [\sigma'] \\
\frac{-(1-p)}{2(M-1)}  & \mbox{when } [\sigma] \downarrow [\sigma'] \\
\frac{1-p}{2(M-1)}  & \mbox{when }  [\sigma] \downarrow -[\sigma'] \\
0  & 
\mbox{otherwise },
\end{cases} \]
The normalized marginal difference of the $p$-lazy Dirichlet walk on the oriented edges of complex as was defined as $\tilde{\mathcal{E}}_n= \left(\frac{M-1}{p(M-2)+1}\right)^n B^n 1_{[e]}$ and they proved \cite{Sayan}: 
\begin{thm}\label{sayan} : 
\begin{itemize}
\item[i.] $B=\frac{p(M-2)+1}{M-1}I- \frac{1-p}{2(M-1)} L_1^{down}$ \\
where $L_1^{down}$ is the Laplacian  with  a constant $1$ weight function on the set of simplexes.
\item[ii.] If $\frac{M-2}{3M-4}<p<1 $, then $\tilde{\mathcal{E}}_{\infty}^{[e_0]}$ exists (for arbitrary starting oriented edge $[e_0]$), the 1-th homology is trivial if and only if for all starting (oriented) simplexes $\tilde{\mathcal{E}}_{\infty} \in \text{im}\delta^*_{1}.  $ Furthermore the dimension of $H_{1}(K)$ equals the dimension of the $$\Span \{\mathbb{P}_{ker \delta_{1}}(\tilde{\mathcal{E}}_{\infty}^ {[e]})\mid e\in S_{1} \}$$ where $\mathbb{P}$ denotes the projection map onto the $\ker \delta_{1}$.
 The same holds for $p=\frac{M-2}{3M-4}$ if the connected $N$-complex is not disorientable and either the dimension of the complex (i.e. $N$) is bigger than one or there are no
disorientable $N$-connected components of constant $(N-1)$-degree.

\item[iii.] If $p\geq \frac{1}{2} $ we have
    \begin{equation}
\dis(\tilde{\mathcal{E}}_t^{[e_0]}, \tilde{\mathcal{E}}_{\infty}^{[e_0]})= O\left((1-\frac{1-p}{2(p(M-2)+1)} \lambda_1\right)^t. 
\end{equation}
     where $ \lambda_1$ is the spectral gap of the $L_1^{down}$.

\end{itemize}
\end{thm}

\begin{rem}
Originally,  the above theorem is expressed for arbitrary $d$-dimension $(1 \leq d \leq N)$ based on the walks among $d$-dimensional "oriented" simplexes by going trough their shared $(d-1)$-simplexes as mentioned before. For general $d$, the formulas are similar as we just need to change all the "1"s to "d" in the above theorem and the coefficients of $L_d^{down}$ and $\lambda_d$ in (i) and (iii) are $\frac{1-p}{(d+1)(M-1)}$ and $\frac{1-p}{(d+1)(p(M-2)+1)}$ respectively.  
\end{rem}
This walk and the up-walk share the same challenges as one has to  characterize the limiting behavior of the difference between two probability measures and not a single one and therefore the limiting behaviour of the walk is dependent on the starting point and irreducibility is unclear.

\paragraph{Normalized Hodge Laplacian}
In another approach, presented in  \cite{Schaub_2020}, the authors introduce a  normalization for the discrete (full) $1$-Hodge Laplacian  on simplicial complexes as a generalisation for the graph case and relate that to a random walk on oriented edges via both upper adjacency and lower adjacency relations. Although they focus on the space of $1$-chains on simplicial complexes, as a natural way of modeling
flows, their idea can be extended to random walks in any dimension. However, in their results the connection  between the rate of convergence of these $1$-walks with the spectrum of the Laplacian (and therefore $1$-th homology) is not explored.  



\subsection{The perspective of signed graphs}
So far, we have discussed the  topologically relevant existing approaches for random walks on higher dimensions and their challenges due to orientation. However, there remain some important questions:
\begin{itemize}
    \item Is it possible to overcome the issue of orientations and the complications that we
already discussed in both of the above random walks?
\item Can we define random walks on higher dimensional simplexes  in the same way as random walks on graphs where one considers
the limiting behavior of a single probability measure, while keeping all the connections to the spectrum
of the Laplacian and homology?
\item When does the limit of the walk become independent of the starting points?
Can we get a clear idea about when the walk on higher simplices is irreducible?
\end{itemize} 
To the best of our knowledge, these questions have not been clearly answered.

The crux of the challenge in going from random walks on graphs to simplicial complexes is orientation. If we can
somehow ignore  orientations  on higher dimensional simplexes and define random walks on the $d$-simplexes and not the oriented $d$-simplicies, while keeping the desired geometric and topological properties of the walk then we return to the graph setting.

Our main contribution  and goal in this paper is to address these questions and state clearly if and when graph-type random walks on simplicial complexes can give rise to the same connection that appears in the theorems above. 
We use signed graphs as well as the Perron-Frobenius theorem \cite{Perron,Frobenius} to achieve
our goal.

There are two main steps to clarify the path to our goal:
\begin{itemize}
\item We want to connect the $d$-th Laplacian spectrum of a $N$-simplicial complex to its $d$-th homology via a random walk on the oriented $d$-simplexes such that the limit of the walk is related to the up/down spectral gap in dimension $d$. 
\item The second and more important aspect is to see how and if we can choose an orientation on the $d$-simplexes such that our proposed walk, becomes independent of the starting point and is as simple and informative as walks on the vertices in graphs. It is important to note that for such a simplified connection we should be sure that for each $d$-simplex just one of its orientations is enough to be a state of the walk in the state space of Markov chain. If such an orientation exists and each $d$-simplex is presented once in the state space, the state space would have the same number of elements as the size of $S_d$. But how do we know that such orientation might exist?
\end{itemize}
  
Having the graph case in mind, to connect the $\Delta_{d}^{up}$ with a random walk on the oriented $d$-simplexes, the desired result would be a (square) stochastic matrix $M$, such that
\begin{equation}\label{Fro} 
    M=\alpha I - \beta \Delta_{d}^{up} 
\end{equation}
where $0 \leq \alpha, \beta \leq1$. For $M$ to be stochastic, first we should make sure that it is non-negative.  Based on such a desired connection between $M$ and $\Delta_{d}^{up}$,  we should look for a choice of orientation on the $d$-simplexes such that in the up-Laplacian matrix all off-diagonal elements are non-positive (otherwise $M$ is not necessarily non-negative). After choosing such an orientation, we can get the desired connection between the limit of the walk which is presented by $M$ and the kernel of the Laplacian which corresponds to its non-trivial homology (due to the Eckmann decomposition).   
 
We know from linear algebra \cite{linear}: 
\begin{itemize}
\item If an irreducible non-negative Matrix $A$ has a positive trace, i.e has at least one element $a_{ij} \succ 0$, then  $A$ is primitive. A primitive matrix is a square non-negative matrix some power of which is positive. Also for finite irreducible Markov chains being primitive is the same as being aperiodic. 

\item  If an irreducible non-negative Matrix $A$ has spectral radius $r$, then $A$ is primitive if and only if:
\begin{equation}   
 \lim_{n\to\infty} (A/r)^n = \frac{v'v^T}{v^{T}v'} 
\end{equation}
Here $v, v'$ are two vectors with strictly positive entries such that $v^TA=rv^T$ and $Av'=rv'$ So $A$ has strictly positive left and right eigenvectors belonging to $r$.
This means that this limit equals to the projection onto the eigenspace corresponding to the spectral radius $r$ This projection is called the Perron projection. 

We note that the above statement for the positive matrices is a result of the Perron–Frobenius theorem which can be generalized to non-negative irreducible matrices with some specific extra assumptions such as primitivity, which is equivalent to aperiodicity.
\end{itemize}
The goal is to connect the limit of the processes to the eigenvalue zero of $\Delta_{d}^{up}$, the minimum eigenvalue of this operator. Therefore with the help of the above theorems we should be sure that the spectral radius of $M$ corresponds to the maximum eigenvalue of $M$. We should have
\begin{equation}
\mid \alpha - \beta \times (d+2) \mid < \alpha
\longleftrightarrow 2\alpha > \beta \times (d+2)
\end{equation}
where $(d+2)$ is the maximum possible eigenvalue of $\Delta_{d}^{up}$ \cite{Jost}.  The left hand side of the arrow follows from the fact that the spectrum of $\Delta_{d}^{up}$  is  $[0, (d+2)]$. Therefore for any $\alpha$ and $\beta$ that satisfy this inequality, the limiting behavior of $M$ would be connected to the maximum eigenvalue of $M$ which gives us the kernel of the $\Delta_{d}^{up}$. On the other side, we can explicitly compute $\alpha$ and $\beta$ based on the definitions of $M$ and $\Delta_{d}^{up}$.
  \begin{rem}
  We note that:
  \begin{itemize}
    \item  The above formula can be generalized to any transition, identity and Laplace operators $M$, $I$ and $L$,  as long as 1) we know the upper bound for $L$ and 2) we are sure that an orientation on the $d$-simplexes exists in such a way that in the Laplacian matrix all off-diagonal elements are non -positive (otherwise $M$ is not necessarily non-negative). In that case $(d+2)$ in the above formula is substituted by the upper bound of $L$. Specifically, the exact same argument works for the down normalized Laplacian which will be discussed later. 
    \item The type of random walk that is used in \cite{Ros} ensures that the off-diagonal elements of $L^{up}$ become non-positive. We recall that in their construction the evolution of the (normalized) expectation process, which is the difference of two probability measures,  is determined by the transition operator $A$ and the state space of this walk is the set of all oriented simplicies but it is not clear if we need all of them to get to  the stationary distribution and the consequent properties. Therefore the above formula (\ref{Fro}) is satisfied (when substituting the corresponding upper bounds for the Laplacian) even though they have not directly used the  Perron–Frobenius theorems to reach theorem \ref{Ros}.
    
    \item The situation for the down walk in \cite{Sayan} is different since the propagation matrix $B$ which (together with the down $d$-Laplacian) satisfies the above formula \ref{Fro}, as shown in theorem \ref{sayan}, might have negative elements and therefore the Perron–Frobenius theorem can not be used directly. However, to show the limiting behaviour of $B$, they use the $T$ operator,  which enforces anti-symmetry between the opposite orientations of $d$-simplexes
    and they prove that for every $t$, $TP^t=B^tT$.
    
    \end{itemize}
     \end{rem}
Now coming back to our main question: is it possible to assign an orientation on the simplexes such that 1) all off-diagonal elements of the Laplacian matrix become non-positive and 2) we need just one of the two orientations of each $d$-simplex
so the walk can be restricted to  half of the oriented $d$-simplexes. If the answer is positive, it is the same as getting rid of orientations, as one of the two orientations in each step of the walk would be sufficient to be used in the matrices in the above formula, in the same manner that we have in graphs. We answer this questions with the help of signed graphs of simplicial complexes in  both up and down walks. 

\textit{Signed graphs:}
A signed graph $G$ is a graph $G$ with positive and negative signs ($+/-$) assigned to its edges. We may switch signs by changing the signs of all edges that are connected to an arbitrarily chosen vertex. A signed graph is called balanced if by switching some vertices, we can make all signs positive. It is antibalanced
if we can make all of the signs negative. 
Signed graphs are used in many areas from modeling biological networks to social relations and signed networks \cite{Harary,neuro,signedc, Signed}. The spectral theory for signed graphs has led to a number of breakthroughs in combinatorial geometry and computer science.  \cite{Computer,mathbreak}
The Laplacian of the signed graph $(G, S)$ is
\begin{equation}
 L^sf(v)= f(v)- \frac{1}{\deg v}\sum_{v'\sim v}^{} s(vv') f(v') 
 \end{equation}
 where $v \sim v'$ means $v$ and $v'$  are connected by an edge and its sign is determined by $s(vv')$.\\
The eigenvalues of $L^s$ are real and lie in the interval $[0, 2]$. The smallest eigenvalue is zero if and only if the signed graph is balanced, and positive otherwise. The largest eigenvalue is $=2$ if and only if the graph is antibalanced\cite{signedL, Jost }.

\paragraph{Signed graph of an oriented simplicial complex}
We can connect the Laplacian of a simplicial complex to the Laplacian of a signed graph and relate their spectrum \cite{Jost}. 
There are different ways to construct signed graphs for an "oriented" simplicial complex based on two factors. The first is the way we define the edges; either to represent the up-adjacency connections of the $d$-simplexes or to represent the down-adjacency connections. The second is how we assign the signs of the edges.

 We consider two types of signed graphs:
\begin{itemize}
 \item The up-signed graph $G_{up}(S_d)$ of the oriented simplicial complex in dimension $d$ is constructed as follows. We take a simplicial N-complex. We choose a $d$ and fix an orientation on $d$ and $(d+1)$- simplexes (we can simply ignore the simplexes of dimensions higher than $(d+1)$). The vertices of this graph are the $d$-oriented simplexes (as already mentioned orientations are arbitrarily chosen but fixed, namely the number of vertices is as many as the elements of $S_d$) and we draw an edge between these oriented $d$-simplexes $[\sigma]$ and $[\sigma']$ whenever $\sigma$ and $\sigma'$ share a coface $\tau$, $ \tau \in S_{d+1}$ such that  $\sigma, \sigma' \subset \tau.$ The sign of an edge is defined as
 \begin{center}
     $- \sgn([\sigma], \partial[\tau])\times \sgn([\sigma'], \partial[\tau])$
     \end{center}
     where $\sgn([\sigma], \partial[\tau])$ is positive if the induced orientation of $\tau$ on $\sigma$, by the boundary operator $\partial,$ is the same as the orientation assigned to $\sigma$ and negative otherwise. Similarly for  $\sgn([\sigma'], \partial[\tau])$.
 \item The down-signed graph $G_{down}(S_d)$: similar to the previous case we choose and fix a $d$ and we assign arbitrary orientations to $d$ and $(d-1)$- simplexes (we simply ignore simplicies of dimensions higher than $d$). The vertices of this graph are the $d$-oriented simplexes and an edge is added whenever simplexes $[\sigma]$ and $[\sigma']$  share a $(d-1)$- dimensional face $\rho$. The sign of an edge is
 \begin{center}
 $
 - \sgn ([\rho], \partial[\sigma])\times \sgn ([\rho], \partial[\sigma'])$
  \end{center}
\end{itemize} 
Now we state our simple but fundamental observation on the following results which answers our main question for the up-walk:  

As shown in \cite{Jost}, It is not possible to orient a simplicial complex in such a way that its up-signed graph becomes balanced. Since for a fix dimension $d$, Laplacian of signed graph of a simplicial complex and the (normalized) up-Laplacian   have direct  positive correlation, as shown in the following formula \cite{Jost}, we can never orient a simplicial complex in such a way that  the signed graph Laplacian matrix get non-positive diagonal elements.  ($ 0 \leq d \leq N-1$).
\begin{equation}
 \Delta_d^{up}= (d+1) L^{s}- d I   
\end{equation}
Equivalently, we can never orient a simplicial complex in such a way that its signed graph becomes a simple graph (without negative edges). Therefore we can never have a graph-type up-random walk on the $d$-simplexes that is sensitive to topology of the complex. And we should always consider both orientations of the $d$-simplicies as the states of the walk to infer topological information.
\begin{rem}
    It worth mentioning that there is a notion of anti-signed graph of a simplicial complex  where its only difference with the above construction for both up and down cases is that the  above formulas of the sign assignment are multiplied by a negative sign. For the up Laplacian, it is straightforward to see that the anti-signed graph of a simplicial complex will be balanced (i.e. becomes the up dual graph) if and only if  the simplicial complex is disorientable. However, if we define a random walk on this dual graph it  corresponds to the disorientability (i.e. the maximum eigenvalue of the up Laplacian) and will not be relevant to the topology of the complex (corresponding to the minimum eigenvalue of the up Laplacian). However, such random walks might be useful for obtaining some Cheeger-type inequalities   from above for simplicial complexes and needs further research.    
\end{rem}
 It remains to explore the down walks and see if a graph-type down walk can tell us about the topology of the complex? 
 To answer this question we again use the signed graphs 
and  have the following observation: 
  \begin{lemma}
     The down $d$-signed graph of a (down-connected) simplicial complex is balanced (i.e. is the same as the down d-dual graph) if and only if the simplicial complex is orientable. Note that this requires $d$ to be the maximum dimension. 
 \end{lemma}
 \begin{proof}
     If the signed graph is balanced, it means that by switching the signs of some of the vertices, all signs of the graph become positive. This signs correspond to a selection of orientations for the $N$ and $(N-1)$-simplexes where the $N$-simplexes induce opposite orientation on their common $(N-1)$-faces. By definition the complex is orientable.  
 \end{proof}

 We note that this can be considered quite restrictive for general simplicial complexes as for orientability the maximum degree of $(N-1)$- simplexes is at most two (i.e. the simplicial complex is non-branching). However, it can be quite useful when it comes to triangulations of orientable manifolds. 

 Assuming orientability, to obtain simple formulas for the down normalized Laplacian of the simplicial complex, we restrict to the case that the degree of all ($N-1$)-simplexes is two. After assuming that all the $N$-simplexes are oriented based on the one used to show that the simplicial complex is orientable we have:
 \begin{equation}
\Delta_{N}^{down}f[\sigma]:= \frac{(N+1)f([\sigma])}{2} -\frac{1}{2}\sum_{\sigma', \exists \rho_{N-1}: \rho \in \sigma \cap \sigma'\ } f([\sigma'])
\end{equation}
With the same assumptions as before we have:
\begin{equation}\label{down}
 \Delta_N^{down}= \frac{N+1}{2} L^{s}   
\end{equation}

By the above formula, since all of the $N$-simplexes are oriented in a compatible manner (i.e. induce opposite orientations on the common faces)  and this is global,  we can think of it as ignoring orientations. If the simplicial complex is orientable, fixing such an orientation on the maximum dimensional simplexes  gets us non-positive off diagonal elements in the Laplacian matrix and a graph-like random walk on $N$-simplexes that is fully informative about the topology of the complex, based on the Perron-Frobenius theorem.
 \begin{rem}
      Our weight function here is different from what is used in \cite{Sayan}  and therefore instead of their upper bound ( $(N+1) \times$ maximum degree of $(N-1)$-simplexes ), the upper bound of $\Delta_{N}^{down}$ will be $(N+1)$ \cite{Horak}. But  we can use any other positive weight function;  however the common idea is that for all of the weight functions, to achieve the upper bound of the $N$-down Laplacian we require that the simplicial complex is disorientable, i.e. we get the upper bound for the corresponding $(N-1)$-up Laplacian.
 \end{rem}

 Now we define a Markov chain $M$ on the $N$-simplexes of the orientable simplicial complex as 
\begin{equation}
 P_{\sigma\rightarrow \sigma'}:=
    \begin{cases}
        p & \text{if } \sigma= \sigma'\\
        \frac{1-p}{\theta} & \text{if } \exists \rho_{N-1} \in \sigma \cap \sigma'\\
        0 & \text{otherwise} 
    \end{cases}
\end{equation}
where $\theta$ is number of non-free faces of $\sigma$. Recall that $\rho$ is a free face for $\sigma$ if it has no other coface except $\sigma$. Since we assume orientability here, each $(N-1)$-face has degree at most two. If the simplicial complex has no free $(N-1)$-face $\theta$ is constant for each $\sigma$ and is equal to $N+1$. For instance, triangulations of closed orientable manifolds have this property. 

\begin{figure}[ht]
\centering
\includegraphics[width=0.35\textwidth]{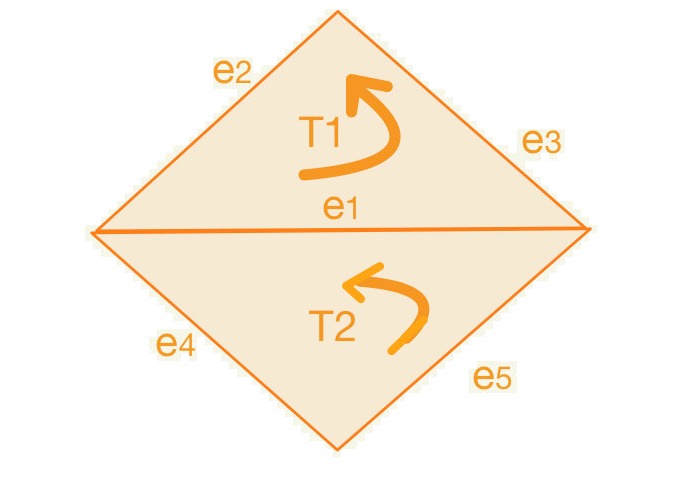}
\caption{Graph-type down random walk on triangles of an orientable 2-d simplicial complex. 
}
\end{figure}
 Here, we start from oriented triangle $T_1$ and with probability $p$, we stay at $T_1$ and with probability ${1-p}$, we jump into oriented $T_2$. We note that the induced orientations of $T_1$ and $T_2$ on their common face $e_1$, is opposite of each other.

For simplicity, we consider the simplicial complex has no free-face, the degree of all of the $(N-1)$-simplexes is two. However, the following results can be extended to general orientable complexes with non-zero boundary. For that we extend the $N$-complex by adding $N$-simplexes to all of the free $(N-1)$-faces (i. e., boundary simplexes) such that the degree of  all of these faces become $2$. Then we use the random walk as described above on this extended N-complex. This is the same idea as the graph random walk with Dirichlet boundary  conditions (meaning
the walk ends  whenever the walker reaches one of the boundary vertices)  that was defined by Chung \cite{Chung} (also see \cite{Sayan}).

\begin{rem}
  Looking closer at the down  random walk presented in \cite{Sayan}, we notice that their random walk is based on a local orientability condition  of the lower adjacent oriented simplexes. Therefore graph-type down random walks on simplicial complexes will be topologically relevant if and only if this local orientability becomes global. 
\end{rem}

By computations from formula 3.12 and 3.14 we have
\begin{equation}
    M= I - \frac{2(1-p)}{N+1} \Delta_{N}^{down} 
\end{equation}

\begin{lemma}
     The spectrum of $M$ is $[2p-1, 1]$ and its minimum is achieved if and only if the $N$-complex has a disorientable $N$-connected component of constant $(N-1)$-degree and
the max is achieved by the cochains in the kernel of $\delta^*_{N-1}$ (i. e. cycles).
\end{lemma}
\begin{proof}
    The first part of the claim can be simply obtained by the formula 3.15 and considering the fact that the spectrum of $\Delta_{N}^{down}$ is contained in $[0, (N+1)]$ where 0 is achieved by the elements in the $ \ker \delta^*_{N-1}$ and the upper bound is achieved if and only if the $N$-complex has a disorientable component of constant $(N-1)$-degree; this can be proven in the same manner as the proof of the lemma 2.5 in \cite{Sayan} and we omit it here.  
\end{proof}

 
  For each $t$, let $p_t^ {\sigma_{0}}= M^t p_{0}^{\sigma_{0}}$ and $p^{\sigma_{0}}_t {(\sigma')}$ is the probability of finding the random walker at $\sigma'$ at time $t$, starting from $\sigma_0$.
\begin{thm}
  The above definition satisfies the following: 
\begin{itemize}
\item[i.]  If $0<p<1 $, then $p_{\infty} = \lim_ {t \rightarrow\infty} p_t$ exists, is independent of the starting point and is equal  to the projection onto the kernel of $ \triangle_N^{down}(= \ker \delta^*_{N-1}$, denoted by $\mathbb{P}_{\ker \delta^*
_{N-1}})$. Also the dimension of $H_{N}(K)$ equals
the dimension of 
 the $\ker \delta^*_{N-1}$ 
 since in the maximum dimension  the $\ker \triangle_N= \ker \triangle_N^{down}$ (as the $\triangle_N^{up}=0$) (compare this theorem with theorem \ref{sayan} which was presented in \cite{Sayan} for any $d$ in  general simplicial complexes). 

 The same holds for $p=0$ if there are no disorientable N-connected component of constant $(N-1)$-degree.

\item[ii.] If $p\geq \frac{1}{2} $ we have
    \begin{equation}
    \dis(p_t, p_{\infty)} = O\left(1-\frac{2(1-p)}{(N+1)} \lambda_N\right)^t.
     \end{equation}
     where $ \lambda_N$ is the spectral gap of the $\triangle_N^{down}$.

\end{itemize}
\end{thm}

\begin{proof}
    To prove the first part, we use again our arguments on the conditions for using the  Perron–Frobenius theorem that was discussed in the previous part. Here, after substituting $\triangle_{N-1}^{up}$ by $\triangle_N^{down}$ in the formula 3.7, we see that:
    
    $2\alpha \geq \beta (N+1) \leftrightarrow 0\leq p<1$.
    
    Also $p=0$ is not problematic except if the connected $N$-complex is disorientable and of constant $(N-1)$-degree. 
    
    For ii) we note that if $p\geq \frac{1}{2}$, $M$ is positive semidefinite and we have:
    
$$\parallel M^t- \mathbb{P}_{\ker \delta
_{N-1}}\parallel= \left \| (I-\frac{2(1-p)}{N+1}\triangle_N^{down})^t\mid_{\im\triangle_N^{down}} \right \|= \left(1-\frac{2(1-p)}{(N+1)} \lambda_N\right)^t$$
    
\end{proof}
\begin{rem}
    Since based on the equation \ref{down}, the down signed-Laplacian and the down normalized Laplacian have direct positive correlation, we can define our random walk  on the down dual graph, exactly in the same way as described above. This walk will be fully informative about the topology of the complex as similar to the normalized graph case in \ref{graph}, we would have $M=I-(1-p)L^s$. Therefore the convergence rate of this walk is determined by the spectral gap of $L^s$ which is directly related to the spectral gap of $\triangle_N^{down}$ according to \ref{down}. However, we used the direct approach to  emphasize the differences with the one that was developed in \cite{Sayan}   
\end{rem}

\section{Discussion}
With the help of the signed graphs of simplicial complexes and the Perron-Frobenius theorem, we can state:   
 
 Topologically relevant up-random walks on simplicial complexes will never be irreducible while down random walks on the simplicial complexes are irreducible if and only if the  walk is happening on  the maximum dimensional simplexes and the complex is orientable.
Therefore the up-walks  can never be simplified (to graph-like walks, reflecting the non-trivial homology) compared to the down walks which can be  simplified when the simplicial complex is orientable. We recall that based on the results in \cite{Sayan}, if we exclude such assumption, we will not have irreducibility of the proposed Markov chain. 
There remain many open questions about random walks on simplicial
complexes and their connections to spectrum of higher order Laplacians and the topology and geometry of the complex. Possible future directions of research include:
\begin{itemize}
    \item What is the continuum limit of these walks on manifolds and can they help us to explore the spectrum of the Hodge Laplacian on manifolds and furthermore to obtain Cheeger-type inequalities in higher dimensions? 
    \item How can we extend some geometric notions such as Ollivier-Ricci from graphs to the signed graphs of general simplicial complexes? Notions that are based on random walks and can be connected to the Laplacian spectrum. 
    \item What are applications of these walks in both irreducible and non-irreducible cases? How much we can extend the applications of graph-walks to the walks in simplicial complexes?
\end{itemize}
\textbf{Acknowledgements}\\
Marzieh Eidi would like to thank Dr. Zahra Eidi for the always insightful talks and her kind advice and support. Also she thanks Prof. Dr.  J\"urgen Jost, director of the Max Planck Institute for Mathematics in the Sciences in Leipzig, for introducing her the notion of the signed graphs which inspired some of ideas in this paper and for the nice memorable talks all these years. Sayan
Mukherjee would like to acknowledge partial funding from NSF DMS 17-13012, NSF 30 BCS 1552848, NSF DBI 1661386, NSF IIS 15-46331, NSF DMS 16-13261, as well as the high-performance
computing partially supported by grant 2016-IDG-1013 from the North Carolina Biotechnology Center as well as the Alexander von Humboldt Foundation, the BMBF and the Saxony State Ministry for
Science.  


\bibliographystyle{unsrt}  
\bibliography{references}  

\end{document}